\documentclass[12pt, reqno]{amsart}
\usepackage{amsmath, amsthm, amscd, amsfonts, amssymb, graphicx, color}
\usepackage[bookmarksnumbered, colorlinks, plainpages]{hyperref}
\usepackage{float}
\usepackage{graphicx,epstopdf}
\usepackage{bbm}
\usepackage{mathrsfs}
\usepackage{dsfont}
\usepackage{amssymb}
\usepackage{csquotes}
\usepackage[dvips]{epsfig}
\usepackage{latexsym}
\usepackage{exscale}
\usepackage[latin1]{inputenc}
\usepackage{setspace}
\newtheorem{theorem}{Theorem}[section]
\newtheorem{lemma}[theorem]{Lemma}

\theoremstyle{definition}

\newtheorem{example}[theorem]{Example}

\theoremstyle{remark}
\newtheorem{remark}[theorem]{Remark}
\numberwithin{equation}{section}
\begin{document}

\setcounter{page}{1}

\title[Bernstein Lethargy Theorem]{Constructing an Element of a Banach Space with Given Deviation from its Nested Subspaces}

\author[A. AKSOY, Q. Peng]{Asuman G\"{u}ven AKSOY$^1$$^*$, Qidi Peng$^2$}

\address{$^{1}$Department of Mathematics, Claremont McKenna College, 850 Columbia Avenue, Claremont, CA  91711, USA.}
\email{aaksoy@cmc.edu}

\address{$^2$Institute of Mathematical Sciences, Claremont Graduate University, 710 N. College Avenue, Claremont, CA  91711, USA.}
\email{qidi.peng@cgu.edu}




\date{Received: xxxxxx; Revised: yyyyyy; Accepted: zzzzzz.
\newline \indent $^{*}$ Corresponding author}

\begin{abstract}
This paper contains two improvements on a theorem of S. N. Bernstein for Banach spaces.  We show that  if $X$ is an arbitrary infinite-dimensional Banach space,  $\{Y_n\}$ is a sequence of strictly nested subspaces of $ X$  and if   $\{d_n\}$ is a non-increasing sequence of non-negative numbers tending to 0, then for any $c\in(0,1]$
  we can find  $x_{c} \in X$,    such that the distance $\rho(x_{c}, Y_n)$ from $x_{c}$ to $Y_n$ satisfies
$$
c d_n \leq \rho(x_{c},Y_n) \leq 4c d_n,~\mbox{for all $n\in\mathbb N$}.
$$
We prove the above inequality by first improving Borodin (2006)'s result for Banach spaces  by weakening his condition on the sequence $\{d_n\}$.  The weakened condition on $d_n$ requires refinement of Borodin's construction to extract  an element in $X$, whose distances from the nested subspaces are  precisely the given values $d_n$.
\end{abstract} \maketitle

\section{Introduction}
For a subspace $A$ of a normed linear space $(X,\|\cdot\|)$,  we define the distance from $x\in X$ to $A$ by
$$\rho(x, A):= \inf \{ \|x-a\|: \,\,a\in A\}.$$
In $1938$, S. N. Bernstein \cite{Bernstein} proved that if $\{d_n \}_{n\ge1}$ is a non-increasing null sequence (i.e., $\lim\limits_{n\to\infty}d_n=0$) of positive numbers, and $\Pi_n$ is the vector space of all real polynomials of degree at most equal to $n$, then there exists a function $f\in C[0,1]$ such that
$$\rho(f, \Pi_n)= d_n,~\mbox{for all $n \ge1$}.$$
This remarkable result is called Bernstein's Lethargy Theorem (BLT) and is used in the constructive theory of functions \cite{Sin}. Then it has been applied to the theory of quasi analytic functions in several complex variables \cite{Ple2}.
Note that the density of polynomials in $C[0,1]$ (the Weierstrass Approximation Theorem) implies that $$\displaystyle \lim_{n \to \infty} \rho(f, \Pi_n) = 0.$$ However, the Weierstrass Approximation Theorem gives no information about the speed of convergence for $\rho(f, \Pi_n)$.  Following the proof of Bernstein \cite{Bernstein}, Timan \cite{Tim} extended his result to an arbitrary system of strictly nested \textit{finite-dimensional} subspaces $\{Y_n\}$.  Later Shapiro \cite{Sha}, replacing $C[0,1]$ with an arbitrary Banach space $(X, \|\cdot\|)$ and $\{\Pi_n\}$  the sequence of $n$-dimensional subspaces of polynomials of degree up to $n$, with a sequence $\{Y_n\}$ where $Y_1 \subset Y_2 \subset \cdots$ are strictly nested \textit{closed subspaces} of $X$, showed that in this setting, for each null sequence $\{d_n\}$ of non-negative numbers, there exists a vector $x\in X$ such that
\begin{equation*}
\label{Shap}
 \rho(x, Y_n) \neq O (d_n),~\mbox{as $n\to\infty$}.
 \end{equation*}
That is, there is no $M>0$ such that
$$\rho(x, Y_n) \leq M d_n,~\mbox{ for all $n\ge1$}.
 $$
 In other words $\rho(x, Y_n)$ decays arbitrarily slowly. This result  was later strengthened by Tyuriemskih \cite{Tyu} who established that the sequence of errors of best approximation from $x$ to $Y_n$, $\{\rho(x,Y_n)\}$, may converge to zero at an arbitrary slow rate up to some choice of $x\in X$. More precisely, for any expanding sequence $\{Y_n\}$ of subspaces of $X$ and for any null sequence $\{d_n\}$ of positive numbers, he constructed an element $x\in X$  such that
 \begin{equation*}
 \label{tyu}
 \lim_{n \rightarrow \infty} \rho(x, Y_n) =0,~\mbox{and}~\rho(x, Y_n) \geq d_n~\mbox{for all $n\ge1$}.
   \end{equation*}
   However, it is also possible that the errors of best approximation $\{\rho(x,Y_n)\}$ may converge to zero arbitrarily fast. For example, it is shown by  Theorem 2.2 in \cite{Al-To} that, under some conditions imposed on $\{Y_n\}$ and $\{d_n\}$, for any null sequence $\{c_n\}$ of positive numbers, there exists an element $x\in X$ such that
   \begin{equation*}
   \label{AlmiraJ}
   \lim_{ n \rightarrow \infty} \frac{\rho(x, Y_n)}{d_n} =0~\mbox{ but}~\frac{\rho(x, Y_n)}{d_n}  \neq  O (c_n)~\mbox{as}~n\to\infty.
   \end{equation*}
    We refer the reader to \cite{De-Hu} for an application of Tyuriemskih's Theorem to convergence of sequence of bounded linear operators and to \cite{Al-Oi} for a generalization of Shapiro's Theorem. We also refer to \cite{Ak-Al,Ak-Pe, Ak-Le,Al-To,Lew,Lew1} for other versions of Bernstein's Lethargy Theorem  and  to \cite{Ak-Le2,Alb,Mich,Ple1,Vas} for Bernstein's Lethargy Theorem for Fr\'{e}chet  spaces. Given an arbitrary Banach space $X$, a strictly increasing sequence $\{Y_n \}$ of subspaces of $X$  and a non-increasing null sequence $\{d_n\}$ of non-negative numbers, one can ask the question whether  there exists $x \in X$ such that $\rho(x, Y_n) = d_n $ for each $n$? For a long time no sequence $\{d_n\}$ of this type was known for which such an element $x$ exists for \textit{all} possible Banach spaces $X$. The only known spaces $X$ in which the answer is always ``yes" are the Hilbert spaces (see \cite{Tyu2}). For a general (separable) Banach space $X$, a solution $x$ is known to exist whenever all $Y_n$  are finite-dimensional (see \cite{Tim}).
Moreover, it is known that  if $X$ has  the above property, then it is reflexive (see \cite{Tyu2}).

In this framework we provide two improvements on a theorem of S. N. Bernstein for Banach spaces. First we improve Borodin's  Theorem 1  in   \cite{Borodin1}. We include the statement as Theorem \ref{thm:Borodin} below. Namely, we obtain the same errors of best approximations as in Theorem \ref{thm:Borodin} below but under a weaker condition on the sequence $\{d_n\}$. Then  we use our first improvement to show that,  if $X$ is an arbitrary infinite-dimensional Banach space, and if   $\{d_n\}$ is a decreasing null sequence of non-negative numbers, then under a natural condition on the subspaces $\{Y_n\}$, for any $c\in(0,1]$, there exists $x_{c} \in X$ such that
$$
c d_n \leq \rho(x_{c},Y_n) \leq 4c d_n,~\mbox{for all $n\ge1$}.
$$
\section{Preliminaries}
\label{preli}
Given a Banach space $X$ and its subspaces $Y_1\subset Y_2 \subset \dots \subset Y_n\subset \dots $, it is clear that $$ \rho(f, Y_1) \geq\rho(f,Y_2)\geq \cdots, ~\mbox{for any $f\in X$}$$ and  thus $\{\rho(f,Y_n)\}_{n\ge1}$ form a non-increasing sequence of errors of best approximation from $f$ to $Y_n$, $n\ge1$.  Furthermore we have:
\begin{description}
\item[Property 1] $ \rho(\lambda x, Y_n)= |\lambda| \rho(x, Y_n)$ for any $x\in X$ and $\lambda\in\mathbb R$;
\item[Property 2] $\rho(x+v, Y_n)= \rho(x, \overline Y_n)$ for any $x\in X$ and $v\in \overline Y_n$;
\item[Property 3] $\rho(x_1+x_2, Y_n)\leq \rho(x_1, Y_n)+\rho(x_2, Y_n)$ and consequently   $$\rho(x_1+x_2, Y_n)\geq |\rho(x_1, Y_n)-\rho(x_2, Y_n)| \quad \mbox{ for any}\quad  x_1,x_2\in X.$$
\end{description}
 Note that we also have: $$ |\rho(x_1, Y_n)-\rho(x_2, Y_n)| \leq \|x_1-x_2\|~ \mbox{ for}~  x_1,x_2 \in X,$$
which  implies that  the mapping $X \longrightarrow \mathbb{R}^+$  defined by $ x\longmapsto \rho(x, Y_n)$ is  continuous and thus properties of continuous mappings such as the intermediate value theorem can be used.

Next, we state a basic  BLT result concerning \textit{finite number} of subspaces, for the proof  of the following lemma we refer the reader to Timan's book  \cite{Tim}.
\begin{lemma}
\label{lemma1}
Let $(X,\|\cdot\|)$ be a normed linear space, $Y_1\subset Y_2\subset\ldots\subset Y_n\subset X$ be a finite system of strictly nested subspaces, $d_1>d_2>\ldots>d_n\ge0$ and $z \in X\backslash Y_n$. Then, there is an element $x\in X$ for which $\rho(x,Y_k)=d_k$ $(k=1,\ldots,n)$,  $\|x\|\le d_1+1$, and $x-\lambda z\in Y_n$ for some $\lambda>0$.
\end{lemma}
It is worth noting that, Borodin \cite{Borodin1} proved Lemma \ref{lemma1} where $(X,\|\cdot\|)$ is assumed to be a Banach space, however with the same proof the result still holds for a normed linear space.

  An element $x\in X$ satisfying $\rho(x, Y_n)=d_n$, $n\ge1$ may exist if the sequence $\{d_n\}$ decreases strictly to zero. Borodin  in \cite{Borodin1} uses the above lemma for Banach space to establish the existence of such an element in case of rapidly decreasing sequences; more precisely, in 2006 he proves the following theorem:
\begin{theorem}[Borodin \cite{Borodin1}, Theorem 1]
  \label{thm:Borodin}
  Let $X$ be an arbitrary Banach space (with finite or infinite dimension),  $Y_1 \subset Y_2
  \subset \cdots$ be an arbitrary countable system of strictly nested subspaces
  in $X$, and fix a numerical sequence $\{d_n\}_{n\ge1}$ satisfying: there exists a natural number $n_0\ge1$ such that
  \begin{equation}
  \label{condition_Borodin}
  d_n >
  \sum_{k = n+1}^\infty d_k~\mbox{for all $n \ge n_0$ at which $d_n > 0$.}
  \end{equation}
  Then there is an element $x \in X$ such that
  \begin{equation}
  \label{rhod}
  \rho(x,Y_n) = d_n, ~\mbox{for all $n\ge1$}.
  \end{equation}
\end{theorem}
The condition (\ref{condition_Borodin}) on the sequence $\{d_n\}$ is the key to the derivation of (\ref{rhod}) in Theorem \ref{thm:Borodin}. Based on this result, Konyagin \cite{Konyagin} in $2013$ takes a further step to show that, for a general non-increasing null sequence $\{d_n\}$, the deviation of $x\in X$ from each subspace $Y_n$ can range in some interval depending on $d_n$.
\begin{theorem}[Konyagin \cite{Konyagin}, Theorem 1]
\label{thm:Konyagin}
Let $X$ be a real Banach space, $Y_1 \subset Y_2
  \subset \cdots$  be a sequence of strictly nested closed linear  subspaces of $X$, and $d_1 \geq d_2 \geq \cdots$  be a non-increasing sequence converging to zero, then there exists an element $x\in X$ such that
the distance $\rho(x, Y_n)$ satisfies the inequalities
  \begin{equation}
  \label{bound_Kongyagin}
   d_n \leq \rho(x, Y_n) \leq 8d_n, \,\,\,\mbox{for $n\ge1$}.
   \end{equation}
\end{theorem}
Note that the condition (\ref{condition_Borodin}) is satisfied when $d_n =(2+\epsilon)^{-n} $ for $\epsilon > 0$ arbitrarily small, however it is not satisfied when $d_n= 2^{-n}$. Of course there are two natural questions to ask:
 \begin{description}
   \item[\textbf{Question $1$}] Is the condition (\ref{condition_Borodin}) necessary for the results in Theorem \ref{thm:Borodin} to hold, or does Theorem \ref{thm:Borodin} still hold for the sequence $d_n=2^{-n}$, $n\ge1$?
   \item[\textbf{Question $2$}] Under the same conditions given in Theorem  \ref{thm:Konyagin}, can the lower and upper bounds of $\rho(x, Y_n)$ in (\ref{bound_Kongyagin}) be improved?
 \end{description}

The aim of this paper is to show that the above two questions have affirmative answers. We have weakened the condition (\ref{condition_Borodin}) in Theorem \ref{thm:Borodin} and obtained the same result given by Theorem \ref{thm:Borodin}.  We were also able to improve the bounds given in  (\ref{bound_Kongyagin}), provided some additional weak subspace condition.

Before we proceed with our results, observe that  in Konyagin's paper \cite{Konyagin} it is assumed that $\{Y_n\}$ are closed
and strictly increasing (we will show this assumption can be weakened to $\overline Y_n\subset Y_{n+1}$). In Borodin's paper \cite{Borodin1}, this is not specified, but from the proof of Theorem \ref{thm:Borodin} it is clear that his proof  works only under assumption that $\overline{Y}_n $ is strictly included
in $Y_{n+1}$.  The necessity of this assumption on subspaces is illustrated by the following example.
\begin{example}\label{EXAM}
Let $X = L^{\infty}[0,1]$ and consider $C[0,1] \subset L^{\infty} [0,1]$ and define the subspaces  of $X$ as follows:

\begin{enumerate}
\item $Y_1 =  \mbox{space of all polynomials}$;

\item $Y_{n+1} = $span$[Y_n \cup \{f_{n}\}]$
where $f_{n} \in C[0,1] \setminus Y_n$, for $n\ge1$.
\end{enumerate}
Observe that by the Weierstrass Theorem we have $\overline {Y}_n = C[0,1] $ for all $n \ge1$.
Take any  $ f \in L^{\infty}[0,1]$ and consider the following cases:
\begin{enumerate}
\renewcommand{\labelenumi}{\alph{enumi})}
\item If  $ f \in C[0,1]$ , then
$$
\rho (f,Y_n) = \rho(f,C[0,1]) = 0~\mbox{for all $n\ge1$}.
$$
\item If $ f \in L^{\infty}[0,1] \setminus C[0,1] $, then
$$
\rho(f, Y_n) = \rho(f, C[0,1]) = d > 0~ \mbox{ (independent of $n$)}.
$$
Note that in above, we have used the fact that $\rho(f,Y_n) = \rho(f,\overline{Y}_n)$.
Hence in this case BLT does not hold.
\end{enumerate}
\end{example}
We will assume that the subspaces $\{Y_n\}$ satisfy $\overline Y_n\subset Y_{n+1}$ for $n\ge1$  for the rest of the paper.

\section{Improvement of Borodin's Result}
Our first main result gives a positive answer to Question $1$, by showing that Theorem \ref{thm:Borodin} can be extended by weakening the strict inequality in (\ref{condition_Borodin}) to a non-strict one:
\begin{equation}
 \label{condition}
 d_n\ge\sum_{k=n+1}^\infty d_k,~\mbox{for every} \,\, n \ge n_0.
 \end{equation}
 Clearly the condition  (\ref{condition}) is weaker than (\ref{condition_Borodin}),  but unlike the condition (\ref{condition_Borodin}), (\ref{condition}) is satisfied by the sequences $\{d_n\}_{n\ge1}$ verifying $d_n=\sum\limits_{k=n+1}^\infty d_k$ for all $n\ge n_0$. For a typical example of such sequence one can take $\{d_n\} =\{2^{-n}\}$. As a consequence, the proof of Theorem \ref{Borodin} requires a finer construction of the element $x$ than that of Theorem \ref{thm:Borodin}. Before proving Theorem \ref{Borodin}, we provide the following technical lemmas.
\begin{lemma}
\label{lemma1'}
Let $(X,\|\cdot\|)$ be a normed linear space and let $Q$ be a subspace of $X$ with $\overline Q\subset X$. Pick two elements $x_1,x_2\in X\backslash \overline Q$ such that $x_2\notin \mbox{span} [\{x_1\}\cup Q]$, and  $\delta\geq0$ such that
\begin{equation}
\label{linearx}
\rho(x_2-\delta x_1,Q)\le \rho(x_2-ax_1,Q),~\mbox{for all $a\ge \delta$.}
\end{equation} Then there exists a nonzero linear functional $f:~X\longmapsto \mathbb R$ such that
\begin{equation}
\label{linearf1}
f(q)=0,~\mbox{for $q\in Q$},~\|f\|=\frac{1}{\rho(x_1,Q)},~f(x_1)=1
\end{equation}
and
\begin{equation}
\label{linearf2}
f(x_2)=\delta-\frac{\rho(x_2-\delta x_1,Q)}{\rho(x_1,Q)}.
\end{equation}
\end{lemma}
\textbf{Proof.} First we note that the existence of $\delta$ satisfying the condition (\ref{linearx}) follows from a standard convexity argument. Now,  let $U$ be the linear subspace of $X$ spanned by $\{x_1\}\cup Q$, then any element in $U$ has unique decomposition of the form $q+\alpha x_1$, for some $q\in Q$, $\alpha\in\mathbb R$. From Hahn-Banach Theorem, we have  the linear functional $g:~U\longmapsto \mathbb R$ defined by
\begin{equation}
\label{def:g}
g(q+\alpha x_1)=\alpha,~\mbox{for any $q\in Q$, $\alpha\in\mathbb R$}
\end{equation}
satisfies
\begin{equation}
\label{linearg}
g(q)=0,~\mbox{for all $q\in Q$},~\|g\|=\frac{1}{\rho(x_1,Q)},~\mbox{and}~g(x_1)=1.
\end{equation}
Since $x_2\notin U$, we then let $S$ be the linear subspace of $X$, spanned by $\{x_2\}\cup U$. Define the \emph{sublinear functional} $p$ to be
\begin{equation}
\label{px}
p(x)=\|g\|\|x\|,~\mbox{for $x\in S$}.
\end{equation}
We then can extend $g$ to a linear functional $\tilde g:~S\longmapsto \mathbb R$, by taking
\begin{equation}
\label{def:nu}
\tilde g(q+\alpha x_1+\beta x_2)=\alpha+\beta \nu,~\mbox{for all $q\in Q$ and $\alpha,\beta\in\mathbb R$},
\end{equation}
where the real number $\nu$ is chosen so that
\begin{equation}
\label{fp}
\tilde g(x)\le p(x),~\mbox{for $x\in S$}.
\end{equation}
From Hahn--Banach Theorem (see \cite{Dunford}, Page 63, inequality (i)), any $\nu$ satisfying the inequalities (\ref{inequality1}) below yields (\ref{fp}). Note that
\begin{equation}
\label{inequality1}
\left\{\begin{array}{ll}
&\nu\ge\max\limits_{q\in Q,~\alpha\in \mathbb R}\left\{-p(-q-\alpha x_1-x_2)-g(q+\alpha x_1)\right\};\\
&\nu\le \min\limits_{q'\in Q,~\alpha'\in \mathbb R}\left\{p(q'+\alpha' x_1+x_2)-g(q'+\alpha' x_1)\right\}.
\end{array}\right.
\end{equation}
Next, we show that $\nu$ can take the value $\delta-\frac{\rho(x_2-\delta x_1,Q)}{\rho(x_1,Q)}$. Using (\ref{def:g}), we know that it is enough to prove
\begin{equation}
\label{inequality2}
-p(-q-\alpha x_1-x_2)-\alpha\le \delta-\frac{\rho(x_2-\delta x_1,Q)}{\rho(x_1,Q)}\le p(q'+\alpha' x_1+x_2)-\alpha',
\end{equation}
for all $q,q'\in Q$ and $\alpha,\alpha'\in\mathbb R$. Equivalently, we will show
\begin{equation}
\label{inequality21}
 \delta-\frac{\rho(x_2-\delta x_1,Q)}{\rho(x_1,Q)}\ge-p(-q-\alpha x_1-x_2)-\alpha,~\mbox{for $q\in Q$ and $\alpha\in\mathbb R$}
\end{equation}
and
\begin{equation}
\label{inequality22}
\delta-\frac{\rho(x_2-\delta x_1,Q)}{\rho(x_1,Q)}\le p(q'+\alpha' x_1+x_2)-\alpha',~\mbox{for $q'\in Q$ and $\alpha'\in\mathbb R$}.
\end{equation}
 To show (\ref{inequality21}), we combine (\ref{px}), (\ref{linearg}),  Property 3 of  $\rho( \cdot , Q)$ and the fact that $q\in Q$ to obtain:
      \begin{eqnarray}
     \label{uppermu}
     &&-p(-q-\alpha x_1-x_2)-\alpha=-\frac{\|q+\alpha x_1+x_2\|}{\rho(x_1,Q)}-\alpha\nonumber\\
     &&\le -\frac{\rho(q+\alpha x_1+x_2,Q)}{\rho(x_1,Q)}-\alpha=-\frac{\rho(\alpha x_1+x_2,Q)}{\rho(x_1,Q)}-\alpha\nonumber\\
     &&\le -\frac{||\alpha+\delta|\rho(x_1,Q)-\rho(x_2-\delta x_1,Q)|}{\rho(x_1,Q)}-\alpha
     \end{eqnarray}
     is true  for all $q\in Q$ and $\alpha\in \mathbb R$.
    Now we have two cases to consider for the value of $\alpha$:
     \begin{enumerate}
     \item
      If
     $
     |\alpha+\delta|\rho(x_1,Q)-\rho(x_2-\delta x_1,Q)\le0,
     $
     then (\ref{uppermu}) yields
     \begin{eqnarray}
     \label{upper11}
     -p(-q-\alpha x_1-x_2)-\alpha&\le&\frac{|\alpha+\delta|\rho(x_1,Q)-\rho(x_2-\delta x_1,Q)}{\rho(x_1,Q)}-\alpha\nonumber\\
     &=&|\alpha+\delta|-\alpha-\frac{\rho(x_2-\delta x_1,Q)}{\rho(x_1,Q)}.
     \end{eqnarray}
     \begin{enumerate}
     \renewcommand{\labelenumi}{\alph{enumi})}

     \item
     If $\alpha+\delta>0$, then (\ref{upper11}) becomes
     \begin{equation}
     \label{upper11'}
     -p(-q-\alpha x_1-x_2)-\alpha\le \delta-\frac{\rho(x_2-\delta x_1,Q)}{\rho(x_1,Q)}.
     \end{equation}
     \item
     If $\alpha+\delta\le0$, then by the assumption (\ref{linearx}) we have
     $$
     \rho(\alpha x_1+x_2,Q) \geq \rho(x_2-\delta x_1,Q).
     $$
     This together with (\ref{uppermu}) and the fact that $\alpha\le -\delta$ implies
     \begin{eqnarray}
     \label{upper11''}
     &&-p(-q-\alpha x_1-x_2)-\alpha\le -\frac{\rho(\alpha x_1+x_2,Q)}{\rho(x_1,Q)}-\alpha\nonumber\\
     &&\le -\frac{\rho(x_2-\delta x_1,Q)}{\rho(x_1,Q)}-\alpha\le \delta-\frac{\rho(x_2-\delta x_1,Q)}{\rho(x_1,Q)}.
     \end{eqnarray}
   \end{enumerate}
     \item
     If
     $
     |\alpha+\delta|\rho(x_1,Q)-\rho(x_2-\delta x_1,Q)>0,
     $
     then
    \begin{enumerate}
    \renewcommand{\labelenumi}{\alph{enumi})}
     \item
      If $\alpha+\delta>0$, we solve for $\alpha$ to obtain
     $$
     -\alpha<\delta-\frac{\rho(x_2-\delta x_1,Q)}{\rho(x_1,Q)}.
     $$
     This fact together with  $p(x) \geq 0$  (see (\ref{px})), straightforwardly implies
     \begin{equation}
     \label{upper12}
     -p(-q-\alpha x_1-x_2)-\alpha\le -\alpha<\delta-\frac{\rho(x_2-\delta x_1,Q)}{\rho(x_1,Q)}.
     \end{equation}
     \item If $\alpha+\delta\le0$, then again (\ref{upper11''}) holds.
    \end{enumerate}
     Therefore (\ref{inequality21}) follows from (\ref{upper11'}), (\ref{upper11''}) and (\ref{upper12}).

      To show (\ref{inequality22}),  again we apply (\ref{px}) and Property 3  of $\rho(\cdot , Q)$, to obtain for all $q'\in Q$ and $\alpha'\in \mathbb R$,
     \begin{eqnarray}
     \label{lowermu}
     &&p(q'+\alpha' x_1+x_2)-\alpha'=\frac{\|q'+\alpha' x_1+x_2\|}{\rho(x_1,Q)}-\alpha'\ge\frac{\rho(\alpha' x_1+x_2,Q)}{\rho(x_1,Q)}-\alpha'\nonumber\\
     &&\ge \frac{||\alpha'+\delta|\rho(x_1,Q)-\rho(x_2-\delta x_1,Q)|}{\rho(x_1,Q)}-\alpha'.
     \end{eqnarray}
     Then two different cases follow, according to the values of $\alpha'$:
    \end{enumerate}
    \begin{enumerate}
     \item If
     $
     |\alpha'+\delta|\rho(x_1,Q)-\rho(x_2-\delta x_1,Q)>0,
     $
      (\ref{lowermu}) implies
     \begin{equation}
     \label{lower11}
     p(q'+\alpha' x_1+x_2)-\alpha'\ge\frac{|\alpha'+\delta|\rho(x_1,Q)-\rho(x_2-\delta x_1,Q)}{\rho(x_1,Q)}-\alpha'.
     \end{equation}
     Now we consider the following subcases
     \begin{enumerate}
     \renewcommand{\labelenumi}{\alph{enumi})}
     \item  If $\alpha'+\delta>0$, then (\ref{lower11}) yields
     \begin{equation}
     \label{lower11'}
     p(q'+\alpha' x_1+x_2)-\alpha'\ge\delta-\frac{\rho(x_2-\delta x_1,Q)}{\rho(x_1,Q)}.
     \end{equation}
    \item  If $\alpha'+\delta\le 0$, then by using (\ref{linearx}) we have
      \begin{eqnarray}
     \label{lower11''}
     &&p(q'+\alpha' x_1+x_2)-\alpha'=\frac{\|q'+\alpha' x_1+x_2\|}{\rho(x_1,Q)}-\alpha'\nonumber\\
     &&\ge\frac{\rho(\alpha' x_1+x_2,Q)}{\rho(x_1,Q)}+\delta\ge \delta -\frac{\rho(x_2-\delta x_1,Q)}{\rho(x_1,Q)}.
     \end{eqnarray}
     \end{enumerate}
     \item
      If
     $
     |\alpha'+\delta|\rho(x_1,Q)-\rho(x_2-\delta x_1,Q)<0,
     $
     and $\alpha'+\delta>0$, then its equivalent inequality
     $$
     -\alpha'>\delta-\frac{\rho(x_2-\delta x_1,Q)}{\rho(x_1,Q)}
     $$
     yields
     \begin{equation}
     \label{lower12}
     p(q'+\alpha' x_1+x_2)-\alpha'\ge -\alpha'>\delta-\frac{\rho(x_2-\delta x_1,Q)}{\rho(x_1,Q)}.
     \end{equation}
     If $\alpha'+\delta\le 0$, then the result is the same as in (\ref{lower11''}).
     \end{enumerate}
     Therefore combining (\ref{lower11'}), (\ref{lower11''}) and (\ref{lower12}), we obtain  (\ref{inequality22}).

 Having proved (\ref{inequality21}) and (\ref{inequality22}), we see  that (\ref{inequality2}) holds. Thus we can take $\displaystyle \nu=\delta-\frac{\rho(x_2-\delta x_1,Q)}{\rho(x_1,Q)}$ in (\ref{def:nu}), to obtain (\ref{linearf2}).

 Now, we show $\|\tilde g\|=\|g\|=\displaystyle \frac{1}{\rho(x_1,Q)}$. This is true, because  the fact that $\tilde g(q)=g(q)$ for $q\in Q$ implies $\|\tilde g\|\ge\|g\|$; and the fact that $\tilde g(x)\le p(x)=\|g\|\|x\|$ leads to $\|\tilde g\|\le \|g\|$. It then follows from (\ref{linearg}) that $\|\tilde g\|=\|g\|=\displaystyle \frac{1}{\rho(x_1,Q)}$.

 Finally, from the Hahn--Banach Extension Theorem $\tilde g$ can be further extended to a linear functional $f:~X\longmapsto \mathbb R$, that satisfies $f(x)=\tilde g(x)$ for $x\in S$ and $\|f\|=\|\tilde g\|$.  $\square$
\begin{remark}
It is worth noting that, with a similar idea one can prove the following :
Let $(X,\|\cdot\|)$ be a normed linear space and let $Q$ be a subspace of $X$ with $\overline Q\subset X$. Pick two elements $x_1,x_2\in X\backslash \overline Q$ such that $x_2\notin \mbox{span} [\{x_1\}\cup Q]$, and assume there exists $\delta\ge0$ such that
$$
\rho(x_2+\delta x_1,Q)\le \rho(x_2+ax_1,Q),~\mbox{for all $a\ge \delta$.}
$$ Then there exists a nonzero linear functional $f:~X\longmapsto \mathbb R$ such that
$$
f(q)=0,~\mbox{for $q\in Q$},~\|f\|=\frac{1}{\rho(x_1,Q)},~f(x_1)=1
$$
and
$$
f(x_2)=-\delta+\frac{\rho(x_2+\delta x_1,Q)}{\rho(x_1,Q)}.
$$
\end{remark}
 \begin{lemma}
 \label{lemma2}
 Let $Q_1,Q_2$ be two subspaces of an arbitrary normed linear space $(Q_3,\|\cdot\|)$, such that $\overline Q_k\subset Q_{k+1}$ for $k=1,2$. Let $\{u_m\}_{m\ge1}$ and $\{v_m\}_{m\ge1}$ be two sequences of non-negative numbers, with $u_m>v_m$ for all $m\ge1$. Then there exist a sequence of elements  $\{q_m\}_{m\ge1}\subset Q_{3}$ and a constant $c\ge1$ such that
 \begin{equation}
 \label{q}
 \rho(q_m,Q_1)=u_m,~\rho(q_m,Q_2)=v_m,~\mbox{for all $m\ge1$}
 \end{equation}
 and
 \begin{equation}
 \label{diffq1}
 \|q_m-q_n\|\le c\left(\max\{u_m,u_n\}-\min\{v_m,v_n\}\right),~\mbox{for all $m,n\ge1$}.
 \end{equation}
 \end{lemma}
 \textbf{Proof.} Since $\overline Q_k\subset Q_{k+1}$ for $k=1,2$, then by Lemma \ref{lemma1} there exists an element $z\in Q_3\backslash \overline Q_2$ such that
 \begin{equation}
 \label{zQ}
 \rho(z,Q_1)=2~\mbox{and}~\rho(z,Q_2)=1.
 \end{equation}
 The fact that $\rho(z,Q_2)=1$ implies that, for any $\varepsilon>0$ arbitrarily small, one can find a corresponding  element $w\in\overline Q_2$ such that
 \begin{equation}
 \label{zw1}
 \|z-w\|=1+\varepsilon.
 \end{equation}
 Since by using Property 3 of $\rho(\cdot,Q_1)$ and (\ref{zQ}),
 $$
 \rho(w,Q_1)\ge\rho(z,Q_1)-\rho(z-w,Q_1)\ge \rho(z,Q_1)-\|z-w\|=1-\varepsilon\neq 0,
 $$ we then let
 \begin{equation}
 \label{mumm}
 \delta_{min}=1+\frac{1+\varepsilon-\rho(z-w,Q_1)}{\rho(w,Q_1)}~\mbox{and}~\delta_{max}=\frac{3+\varepsilon}{\rho(w,Q_1)}.
 \end{equation}
 It is clear that $1\le\delta_{min}\le \delta_{max}$, thanks to the fact that $\rho(z-w,Q_1)\le\|z-w\|=1+\varepsilon$ and Property 3 of $\rho(\cdot,Q_1)$. Next observe from (\ref{mumm}) and (\ref{zQ}) that,
 \begin{equation}
 \label{rholess}
 \rho\left(z-\delta_{min}w,Q_1\right)\le \rho(z-w,Q_1)+(\delta_{min}-1)\rho(w,Q_1)=1+\varepsilon
 \end{equation}
 and
 \begin{equation}
 \label{rhogreater}
 \rho\left(z-\delta_{max}w,Q_1\right)\ge \delta_{max}\rho(w,Q_1)-\rho(z,Q_1)=3+\varepsilon-2=1+\varepsilon.
 \end{equation}
 The mapping $\lambda \longmapsto \rho(z-\lambda w,Q_1)$ is continuous, then by (\ref{rholess}), (\ref{rhogreater}) and the intermediate value theorem, there exists a set of $\{\delta_i\}\subset\left[\delta_{min},\delta_{max}\right]$ such that $
 \delta_i<\delta_{i+1}$ for all $i$ and
  \begin{equation}
 \label{zw}
 \rho(z-\delta_i w,Q_1)=1+\varepsilon.
 \end{equation}
 If the set $\{\delta_i\}$ is finite, say $\{\delta_i\}=\{\delta_1,\delta_2,\ldots,\delta_K\}$, then we denote
 \begin{equation}
 \label{delta1}
 \delta=\delta_K.
 \end{equation}
 If the set $\{\delta_i\}$ is an infinite sequence, then since it is strictly increasing and bounded, and the mapping $\lambda\rightarrow\rho(z-\lambda w,Q_1)$ is continuous, there exists a limit $\delta^*=\lim\limits_{i\to\infty}\delta_i$ such that $\delta^*\in\left[\delta_{min},\delta_{max}\right]$ and
$$
 \rho(z-\delta^* w,Q_1)=1+\varepsilon.
$$
In this case we define
\begin{equation}
 \label{delta2}
 \delta=\delta^*.
 \end{equation}
 It follows from (\ref{delta1}) and (\ref{delta2}) that
 \begin{equation}
 \label{bounddelta1}
 1+\varepsilon=\rho(z-\delta w,Q_1)\le \rho(z-a w,Q_1),~\mbox{for all $a\in[\delta,\delta_{max}]$}.
 \end{equation}

 \textbf{Construction of $\{q_m\}_{m\ge1}$:} The fact that $z\in Q_3\setminus \overline{Q}_2$ results in $$z\notin \mbox{span}[\{w\}\cup Q_1] \subseteq Q_2.$$
 Observe the triangle inequality
 $$
 \rho\left(z-aw,Q_1\right)\ge a\rho(w,Q_1)-\rho(z,Q_1)\ge 1+\varepsilon,~\mbox{for all $a\ge \delta_{max}$}.
 $$
 Then from the above inequality and (\ref{bounddelta1}), we know that
 $$
 1+\varepsilon=\rho(z-\delta w,Q_1)\le \rho(z-aw,Q_1),~\mbox{for all $a\ge \delta$}.
 $$
Then we can apply Lemma \ref{lemma1'} to confirm the existence of a nonzero real-valued linear functional $f:~Q_3\rightarrow \mathbb R$ such that
\begin{eqnarray}
\label{f1}
&&Q_1\subset ker f,~\|f\|=\frac{1}{\rho(w,Q_1)},~f(w)=1,\nonumber\\
&&f(z)=\delta-\frac{\rho(z-\delta w,Q_1)}{\rho(w,Q_1)}=\delta-\frac{1}{\rho(w,Q_1)}.
\end{eqnarray}
We define
 \begin{equation}
 \label{x}
 x_1=(f(z)-\delta)w,~x_2=z-f(z)w.
 \end{equation}
 We will show that the sequence $\{q_m\}_{m\ge1}$ satisfying (\ref{q}) and (\ref{diffq1}) can be found in $\mbox{span}[\{x_1,x_2\}]$ (see (\ref{q1}) below). Using (\ref{x}) and (\ref{zQ}) we obtain
 \begin{equation}
 \label{upper1}
 \rho(v_mx_2,Q_2)=v_m.
 \end{equation}
 First, by using (\ref{x}), (\ref{zw}) and the fact that $v_m<u_m$, we have
 \begin{equation}
 \label{upperu}
\rho(v_mx_2+v_mx_1,Q_1)=|v_m|\rho(z- \delta w,Q_1)=v_m\le u_m.
 \end{equation}
 And, since the kernel $ker f$ satisfies
 \begin{equation}
 \label{kerf}
 \rho(x,ker f)=\frac{|f(x)|}{\|f\|},~\mbox{for all $x\in Q_3$},
 \end{equation}
 then by using  (\ref{kerf}),   (\ref{x}) and (\ref{f1}), we obtain
\begin{eqnarray}
\label{loweru}
&&\rho(v_mx_2+u_mx_1,Q_1)\ge \rho(v_mx_2+u_mx_1,ker f)\nonumber\\
&&=\frac{|f(v_mx_2+u_mx_1)|}{\|f\|}=\rho(w,Q_1)\left|v_mf(z-f(z)w)+u_mf((f(z)-\delta w))\right|\nonumber\\
&&=\rho(w,Q_1)u_m\left|f(z)-\delta\right|=\rho(w,Q_1)u_m\left|\delta-\frac{1}{\rho(w,Q_1)}-\delta\right|\nonumber\\
&&= \rho(w,Q_1)u_m\left(\frac{1}{\rho(w,Q_1)}\right)=u_m.
\end{eqnarray}
Since the mapping $\lambda \longmapsto \rho(v_mx_2+\lambda x_1,Q_1)$ is continuous, it follows from (\ref{upperu}), (\ref{loweru})  and the intermediate value theorem that there is a real number $\mu_m\in[v_m,u_m]$ such that
$$
\rho(v_mx_2+\mu_m x_1,Q_1)=u_m~\mbox{and}~\rho(v_mx_2+\mu_m x_1,Q_2)=v_m.
$$
We then denote by
\begin{equation}
\label{q1}
q_m=v_mx_2+\mu_m x_1=v_m(z-f(z)w)+\mu_m(f(z)-\delta)w,~\mbox{for all $m\ge1$}.
\end{equation}
As a consequence (\ref{q}) holds:
$$
\rho(q_{m},Q_1)=u_m,~\rho(q_{m},Q_2)=v_m,~\mbox{for all $m\ge1$}.
$$
Now we show (\ref{diffq1}) holds.  To this end we first state the following 2 evident facts:
\begin{enumerate}
\item If $\mu_m\in[v_m,u_m]$ for any $m\ge1$, then we have
\begin{equation}
\label{diffmu}
|\mu_m-\mu_n|\le \max\{u_m,u_n\}-\min\{v_m,v_n\},~\mbox{for any $m,n\ge1$}.
\end{equation}
\item For any 4 real numbers $u_m,v_m,u_n,v_n$ such that $u_m>v_m,~u_n>v_n$, the following inequality holds:
\begin{equation}
\label{vmvn}
|v_m-v_n|\le \max\{u_m,u_n\}-\min\{v_m,v_n\}.
\end{equation}
\end{enumerate}
It results from (\ref{q1}), the triangle inequality, (\ref{diffmu}) and (\ref{vmvn}) \\that for any $m,n\ge1$,
\begin{eqnarray*}
\label{dq}
\|q_m-q_n\|&\le& |v_m-v_n|\|z-f(z)w\|+|\mu_m-\mu_n||f(z)-\delta|\|w\|\nonumber\\
&\le& c\left(\max\{u_m,u_n\}-\min\{v_m,v_n\}\right),
\end{eqnarray*}
where
$$
c=\max\left\{\|z-f(z)w\|,|f(z)-\delta|\|w\|\right\}\ge |f(z)-\delta|\|w\|\ge\frac{\|w\|}{\rho(w,Q_1)}\ge1,
$$ The last inequality follows from (\ref{f1})  and hence (\ref{diffq1}) holds.
 Lemma \ref{lemma2} is proved by combining (\ref{q}) and (\ref{diffq1}). $\square$


Now, we are ready to prove the following  theorem which  improves the theorem  of Borodin \cite{Borodin1}.
\begin{theorem}
\label{Borodin}
Let $X$ be an arbitrary infinite-dimensional Banach space, $\{Y_n\}_{n\ge1}$ be an arbitrary system of strictly nested subspaces with the property $\overline Y_n \subset Y_{n+1}$ for all $n\ge1$, and let the non-negative numbers $\{d_n\}_{n\ge1}$ satisfy the following property: there is an integer $ n_0 \ge 1 $ such that
$$
 d_n\ge\sum_{k=n+1}^\infty d_k,~\mbox{for every} \,\, n \ge n_0.
$$
 Then there exists an element $x\in X$ such that $\rho(x,Y_n)=d_n$ for all $n\ge1$.
\end{theorem}
\textbf{Proof.} We start by observing the following four cases:
\begin{enumerate}

\item If $n_0\ge2$, the problem is easily converted to the case $n_0=1$: having constructed an element $z$ with $\rho(z,Y_n)=d_n$ for all $n\ge n_0$, we can use Lemma \ref{lemma1} to construct an element $x$ with $\rho(x,Y_k)=d_k$ at $k=1,\ldots,n_0$ and such that $x-\lambda z\in Y_{n_0}$ for some $\lambda>0$. But then observe that
$$
d_{n_0}=\rho(x,Y_{n_0})=\rho(\lambda z,Y_{n_0})=\lambda d_{n_0};
$$
therefore, $\lambda=1$ and
$$\rho(x,Y_n)=\rho(z,Y_n)=d_n~\mbox{ for all $n\ge n_0$}.
$$
Finally $\rho(x,Y_n)=d_n$ for all $n\ge1$.
\item
If $Y_1=\{0\}$, we first convert the problem to the case $n_0=2$ (notice that $Y_2\ne\{0\}$), then by using the above argument, convert the problem to the case $n_0=1$.
\item If $d_n=0$ starting from some $n$, then the desired element exists by applying Lemma \ref{lemma1} to $X=Y_n$ and the subspaces $Y_1\subset Y_2\subset\ldots\subset Y_{n-1}$ within it.
\item Thus we assume $n_0=1$, $Y_1 \neq \{0\}$ and  $d_n>0$, $n\ge1$ for the rest of the proof.
\end{enumerate}
For each $j\ge1$, we define
$$
\tau_{j}=\sum_{k=j+1}^\infty d_{k}.
$$
In view of the above assumptions, we know $\tau_{j}>0$ and the sequence $\{\tau_{j}\}_{j\ge1}$ is monotonically decreasing to $0$. Since $Y_j\ne \{0\}$, then for any integers $j$, $n$ with $1\le j\le n$, we can set
$$
Q_1=\{0\},~Q_2=Y_j,~Q_3=Y_{j+1};~u_n=1+\frac{\tau_n}{2^jd_j},~v_n=1
$$
and
apply Lemma \ref{lemma2} to obtain, the existence of a sequence $\{q_{j,n}\}_{n\ge j}\subset Y_{j+1}\backslash Y_j$ and a constant $c>0$ such that
\begin{equation}
\label{qjnY}
\rho(q_{j,n},\{0\})=1+\frac{\tau_n}{2^jd_j},~\rho(q_{j,n},Y_j)=1
\end{equation}
and
\begin{equation}
\label{diffq}
\|q_{j,m}-q_{j,n}\|\le c\frac{\tau_m}{2^jd_j},~\mbox{for all $n\ge m\ge j$}.
\end{equation}
Now we fix $n\ge1$. Take $\lambda_{n,n}:=d_n$. We see clearly from (\ref{qjnY}) that
\begin{equation}
\label{identqn}
\rho(\lambda_{n,n}q_{n,n},Y_n)=\lambda_{n,n}\rho(q_{n,n},Y_n)=d_n.
\end{equation}
Since $\rho(q_{j,n},Y_j)=1>0$ for $j\le n$ (see (\ref{qjnY})), then there is a nonzero real-valued linear functional (see \cite{Dunford}, Page $64$) $f_{j,n}:~Y_{j+1}\longmapsto\mathbb R$ such that
\begin{equation}
\label{f}
Y_j\subset ker f_{j,n}~\mbox{and}~\|f_{j,n}\|=f_{j,n}(q_{j,n})=1.
\end{equation}
We first assume $f_{n-1,n}(q_{n,n})\ge0$. Then first  by using (\ref{qjnY}) and (\ref{f}), we obtain
\begin{eqnarray}
\label{upperqn}
&&\rho(\lambda_{n,n}q_{n,n},Y_{n-1})\le \rho(\lambda_{n,n}q_{n,n},\{0\})=\lambda_{n,n}\rho(q_{n,n},\{0\})\nonumber\\
&&=d_n\left(1+\frac{\tau_n}{2^nd_n}\right)\le  d_n+\tau_n=\sum_{k=n}^{\infty}d_k\le d_{n-1};
\end{eqnarray}
and then, by the properties of $f_{n-1,n}$ in (\ref{f}), we see
\begin{equation}
\label{lowerqn}
\rho(\lambda_{n,n}q_{n,n}+d_{n-1}q_{n-1,n},Y_{n-1})\ge d_{n-1}f_{n-1,n}(q_{n-1,n})= d_{n-1}.
\end{equation}
It follows from Property 3, (\ref{upperqn}) and (\ref{lowerqn}) that the mapping
$$
\lambda\longmapsto\rho(\lambda_{n,n}q_{n,n}+\lambda q_{n-1,n},Y_{n-1})
$$
 is continuous, through which the image of $[0,d_{n-1}]$ contains $d_{n-1}$. Then by the intermediate value theorem there exists $\lambda_{n-1,n}\in [0,d_{n-1}]$ such that
\begin{equation}
\label{rhodn}
\rho(\lambda_{n,n}q_{n,n}+\lambda_{n-1,n}q_{n-1,n},Y_{n-1})= d_{n-1}.
\end{equation}
If in contrast $f_{n-1,n}(q_{n,n})<0$, we still have (\ref{upperqn}) and moreover,
\begin{eqnarray*}
\rho(\lambda_{n,n}q_{n,n}-d_{n-1}q_{n-1,n},Y_{n-1})&\ge& -\lambda_{n,n}f_{n-1,n}(q_{n,n})+d_{n-1}f_{n-1,n}(q_{n-1,n})\\
&=& d_{n-1}.
\end{eqnarray*}
This implies that we can find  $\lambda_{n-1,n}$ which lies in the interval $[-d_{n-1},0]$ and satisfies the above equation (\ref{rhodn}). Furthermore, in both cases of $f_{n-1,n}(q_{n,n})$, by the fact that $q_{n-1,n}\in Y_n$ and (\ref{identqn}), we have
$$
\rho(\lambda_{n,n}q_{n,n}+\lambda_{n-1,n}q_{n-1,n},Y_{n})= \rho(\lambda_{n,n}q_{n,n},Y_{n})= d_{n}.
$$
For any $2\le k\le n$, assume that we can find the real numbers
\begin{equation}
\label{int-lambda}
\lambda_{n,n}\in[-d_n,d_n],~\lambda_{n-1,n}\in[-d_{n-1},d_{n-1}],\ldots,\lambda_{k,n}\in[-d_k,d_k]
\end{equation}
such that
$$
\rho(\lambda_{n,n}q_{n,n}+\lambda_{n-1,n}q_{n-1,n}+\ldots+\lambda_{k,n}q_{k,n},Y_m)=d_m,~\mbox{for}~m=k,k+1,\ldots,n.
$$
Let
$$
z_{k,n}=\lambda_{n,n}q_{n,n}+\ldots+\lambda_{k,n}q_{k,n}.
$$
Without loss of generality, suppose $f_{k-1,n}(z_{k,n})\ge0$. Then first  by using the triangle inequality, (\ref{int-lambda}) and (\ref{qjnY}), we obtain
\begin{eqnarray}
\label{zY}
&&\rho(z_{k,n},Y_{k-1})\le\|z_{k,n}\|\le\sum_{j=k}^n|\lambda_{j,n}|\|q_{j,n}\|\le \sum_{j=k}^nd_j\left(1+\frac{\tau_{n}}{2^jd_j}\right)\nonumber\\
&&\le \sum_{j=k}^nd_j+\tau_{n}\sum_{j=1}^{\infty}2^{-j}=\sum_{j=k}^nd_j+\tau_n\le d_{k-1}
\end{eqnarray}
and
$$
\rho(z_{k,n}+d_{k-1}q_{k-1,n},Y_{k-1})\ge d_{k-1}f_{k-1,n}(q_{k-1,n})= d_{k-1}.
$$
Therefore, there is $\lambda_{k-1,n}\in [0,d_{k-1}]$ such that
$$
\rho(z_{k,n}+\lambda_{k-1,n}q_{k-1,n},Y_{k-1})=d_{k-1}.
$$
(If $f_{k-1,n}(z_{k,n})<0$, then the number $\lambda_{k-1,n}$ must be found in $[-d_{k-1},0]$.) Furthermore,
$$
\rho(z_{k,n}+\lambda_{k-1,n}q_{k-1,n},Y_{m})=\rho(z_{k,n},Y_m)=d_{m}~\mbox{for}~m=k,\ldots,n.
$$
Continuing this procedure until $k=1$ is included, we obtain the element
$$
x_{n,n}=\lambda_{n,n}q_{n,n}+\ldots+\lambda_{1,n}q_{1,n},
$$
for which $\rho(x_{n,n},Y_k)=d_k$ and $|\lambda_{k,n}|\le d_k$ for $1\le k\le n$. Using the usual diagonalization process, we choose a sequence $\Lambda$ of indices $n$ such that, for all $k\ge1$, $\lambda_{k,n}$ converges to the limit $\lambda_k$ as $n\rightarrow\infty$, $n\in\Lambda$. We then claim that $|\lambda_k|\le d_k$ and the limit $\lim\limits_{n\rightarrow\infty}\sum\limits_{k=1}^n\lambda_kq_{k,n}$ exists in $X$. As a matter of fact, by using the triangle inequality, we have for $m\le n$,
\begin{eqnarray}
\label{dfflambdaq}
&&\Big\|\sum_{k=1}^n\lambda_{k}q_{k,n}-\sum_{k=1}^m\lambda_{k}q_{k,m}\Big\|=\Big\|\sum_{k=1}^m\lambda_{k}(q_{k,n}-q_{k,m})+\sum_{k={m+1}}^n\lambda_{k}q_{k,n}\Big\|\nonumber\\
&&\le \sum_{k=1}^m|\lambda_{k}|\|q_{k,n}-q_{k,m}\|+\sum_{k={m+1}}^n|\lambda_{k}|\|q_{k,n}\|.
\end{eqnarray}
First by using $|\lambda_k|\le d_k$ and (\ref{diffq}), we obtain
\begin{equation}
\label{dfflambdaq1}
\sum_{k=1}^m|\lambda_{k}|\|q_{k,n}-q_{k,m}\|\le c\tau_m\sum_{k=1}^m2^{-k}\xrightarrow[n,m\rightarrow\infty]{}0,
\end{equation}
 next from (\ref{zY}) we see
\begin{equation}
\label{dfflambdaq2}
\sum_{k=m}^n|\lambda_k|\|q_{k,n}\|\le d_{m-1}\xrightarrow[n,m\rightarrow\infty]{}0.
\end{equation}
Combining (\ref{dfflambdaq}), (\ref{dfflambdaq1}) and (\ref{dfflambdaq2}) yields
$$
\Big\|\sum_{k=1}^n\lambda_{k}q_{k,n}-\sum_{k=1}^m\lambda_{k}q_{k,m}\Big\|\xrightarrow[n,m\rightarrow\infty]{}0.
$$
i.e., $\Big\{\sum\limits_{k=1}^n\lambda_kq_{k,n}\Big\}_{n\ge1}$ is a Cauchy sequence in the Banach space $X$, therefore it has a limit in $X$.

Further, we claim that the element
$$
x:=\lim_{n\rightarrow\infty}\sum_{k=1}^n\lambda_kq_{k,n}
$$
is the limit of the sequence $\{x_{n,n}\}_{n\in\Lambda}$ as $n\to\infty$. By using the facts that $\|q_{k,n}\|\le 2$, $|\lambda_{k,n}|\le d_k$ for all $n\ge1$, $k\le n$, and $|\lambda_k|\le d_k$ for all $k\ge1$, we obtain
\begin{eqnarray*}
&&\|x-x_{n,n}\|\le \|x-x_n\|+\|x_n-x_{n,n}\|\\
&&\le \|x-x_n\|+\sum_{k=1}^n|\lambda_{k,n}-\lambda_k|\|q_{k,n}\|\\
&&\le \|x-x_n\|+\sum_{k=1}^{N}|\lambda_{k,n}-\lambda_k|\|q_{k,n}\|+\sum_{k=N+1}^n(|\lambda_{k,n}|+|\lambda_k|)\|q_{k,n}\|\\
&&\le \|x-x_n\|+2\sum_{k=1}^{N}|\lambda_{k,n}-\lambda_k|+4\sum_{k=N+1}^nd_k\\
&&\le \|x-x_n\|+2 N\max_{1\le k\le N,N<n}|\lambda_{k,n}-\lambda_k|+4d_N\\
&&\xrightarrow[n\rightarrow\infty,n\in\Lambda,N\rightarrow\infty]{}0.
\end{eqnarray*}
Finally, by the continuity of $x\longmapsto\rho(x,Y_k)$, we obtain
$$
\rho(x,Y_k)=\lim_{n\rightarrow\infty,n\in\Lambda}\rho(x_{n,n},Y_k)=d_k,~\mbox{for all}~k\ge1.~\square
$$
\begin{remark}
In the above theorem, we improve Borodin's  hypothesis comparing $d_n$ with the tail of the sequence from a strict inequality to $\geq$. This  means that our result  applies when $d_n = \displaystyle \frac{1}{2^n}$, which was not the case in Borodin's formulation. However, our proof is not a mere existence result. We also give an explicit construction  of an element $x\in X$, where the sequence $\rho(x, Y_n)$ is exactly  the sequence $d_n$.  The condition $\overline{ Y_n} \subset Y_{n+1}$ does not come at the expense of our assumption to weaken the condition  on  the sequence $d_n$. This  is a natural condition. To clarify the reason why almost all Lethargy Theorems have this condition on the subspaces, we gave a simple example before Theorem 3.4, (see Example \ref{EXAM}).
\end{remark}
\section{Improvement of Konyagin's result }
As an  application of Theorem \ref{Borodin}, we derive our second main result as an improvement of Konyagin's result in \cite{Konyagin}, which provides a positive answer to Question $2$ that we posed in Section \ref{preli}. In Konyagin \cite{Konyagin} the interval for the error of best approximation  was $[d_n, 8d_n]$, according to our Theorem \ref{Konyagin} below this range can be reduced to $[d_n, 4d_n]$, under a natural subspace condition on $\{Y_n\}$. Notice that in the theorem below, we exclude the cases $(1)$ $d_1\ge d_2\ge\ldots\ge d_{n_0}>d_{n_0+1}=0$ for some $n_0\ge1$; $(2)$ All $\{Y_n\}$ are finite-dimensional or Hilbert spaces; $(3)$ $\{d_n\}$ satisfies the condition (\ref{condition}). It has been proved that in Cases $(2)$ and $(3)$, there is an element $x\in X$ such that $\rho(x,Y_n)=d_n$ for $n\ge1$. In Case $(1)$, analogous to the remark in Konyagin \cite{Konyagin}, Page 206, for every $\varepsilon>0$ arbitrarily small, we can define
$$
d_n'=\left(1+\frac{(n_0-n)\varepsilon}{n_0}\right)d_n,~\mbox{for $n=1,2,\ldots,n_0$}.
$$
We then observe that $$d_1'>d_2'>\ldots>d_{n_0}'=d_{n_0},\quad \mbox{ and}\quad  Y_1\subset Y_2\subset\ldots\subset Y_{n_0}\subset Y_{n_0+1}.$$ Applying Lemma \ref{lemma1} to $\{Y_n\}_{1\le n\le n_0}$ and $\{d_n'\}_{1\le n\le n_0}$, we obtain an element $x\in Y_{n_0+1}$ such that
$$
\rho(x,Y_n)=d_n'\in[d_n,(1+\varepsilon)d_n],~\mbox{for $n=1,2,\ldots,n_0$},
$$
and
$$
\rho(x,Y_n)=0=d_n,~\mbox{for all $n>n_0$}.
$$
Therefore $d_n\le \rho(x,Y_n)\le (1+\varepsilon)d_n$ for all $n\ge1$.
\begin{theorem}
\label{Konyagin}
Let $X$ be an infinite-dimensional Banach space, $\{Y_n\}$ be a system of strictly nested subspaces of $X$ satisfying the condition $\overline{Y}_n \subset Y_{n+1}$ for all $n\ge1$. Let $\{d_n\}_{n\ge1}$ be a non-increasing null sequence of strictly positive numbers. Assume that there exists an extension $\{(\tilde d_n,\tilde Y_n)\}_{n\ge1}\supseteq \{(d_n,Y_n)\}_{n\ge 1}$ satisfying: $\{\tilde d_n\}_{n\ge1}$ is a non-increasing null sequence of strictly positive values, $\overline{\tilde Y}_n\subset \tilde Y_{n+1}$ for $n\ge1$; and there are some integer $i_0\ge1$ and a constant $K>0$ such that
$$
\{K2^{-n}\}_{n\ge i_0}\subseteq \{\tilde d_n\}_{n\ge 1}.
$$ Then for any $c\in(0,1]$, there exists an element $x_c\in X$ (depending on $c$) such that
\begin{equation}
\label{xc}
cd_n\le \rho(x_c,Y_n)\le 4c d_n,~\mbox{for $n\ge1$}.
\end{equation}
\end{theorem}
\textbf{Proof.} We first show (\ref{xc}) holds for $c=1$.

By assumption, there is a subsequence $\{n_i\}_{i\ge i_0}$ of $\mathbb N$ such that
$$
\tilde d_{n_i}=K2^{-i},~\mbox{for $i\ge i_0$}.
$$
Since the sequence $\{\tilde d_n\}_{n=1,2,\ldots,n_{i_0}-1}\cup\{\tilde d_{n_i}\}_{i\ge i_0}$ satisfies the condition (\ref{condition}) and $\overline{\tilde Y}_{n}\subset \tilde Y_{n+1}$ for all $n\ge1$, then we can apply Theorem \ref{Borodin} to get $x\in X$ so that
\begin{equation}
\label{distYd}
\rho(x,\tilde Y_n)=\tilde d_n,~\mbox{for $n=1,\ldots,n_{i_0}-1$},~\mbox{and}~\rho(x,\tilde Y_{n_i})=\tilde d_{n_i},~\mbox{for all $i\ge i_0$}.
\end{equation}
Therefore for any integer $n\ge1$,
\begin{description}
\item[Case 1] if $n\le n_{i_0}-1$ or $n=n_i$ for some $i\ge i_0$, it then follows from (\ref{distYd}) that
$$
\rho(x,\tilde Y_n)= \tilde d_n;
$$
\item[Case 2] if $n_i<n< n_{i+1}$ for some $i\ge i_0$, then the facts that $\{\tilde d_n\}$ is non-increasing and $\tilde Y_{n_i}\subset \tilde Y_{n}\subset \tilde Y_{n_{i+1}}$ lead to
$$
\rho(x,\tilde Y_n)\in\left(\rho(x,\tilde Y_{n_{i+1}}),\rho(x,\tilde Y_{n_i})\right)=\left(K2^{-(i+1)},K2^{-i}\right)
$$
and
$$
\tilde d_n\in\left[K2^{-i},K2^{-i+1}\right].
$$
It follows that
$$
\frac{\rho(x,\tilde Y_n)}{\tilde d_n}\in\left(\frac{K2^{-i-1}}{K2^{-i+1}},\frac{K2^{-i}}{K2^{-i}}\right)=\left(\frac{1}{4},1\right).
$$
\end{description}
Putting together the 2 above cases yields
$$
\frac{1}{4}\tilde d_n\le \rho(x,\tilde Y_n)\le \tilde d_n~\mbox{for all $n\ge1$}.
$$
For $c\in(0,1]$, taking $x_c=4cx$ in the above inequalities, we obtain
$$
c\tilde d_n\le \rho(x_c,\tilde Y_n)\le 4c\tilde d_n,~\mbox{for all $n\ge1$}.
$$
Remembering that $\{(d_n,Y_n)\}_{n\ge1}\subseteq\{(\tilde d_n,\tilde Y_n)\}_{n\ge1}$, we then necessarily have
\begin{equation*}
\label{notequal}
c d_n\le \rho(x_c, Y_n)\le 4c d_n,~\mbox{for all $n\ge1$}.
\end{equation*}
Hence Theorem \ref{Konyagin} is proved. $\square$
\begin{remark}
The subspace condition given in Theorem \ref{Konyagin} states that the nested sequence $\{Y_n\}$ has \enquote{enough gaps} so that the sequence $$\{(d_n',Y_n')\}_{n\ge1}=\{(d_n,Y_n)\}_{n\ge1}\cup\{(K2^{-i},\tilde Y_{n_i})\}_{i\ge i_0}$$ satisfies $d_n'\ge d_{n+1}'\to 0$ and $\overline {Y_n'}\subset Y_{n+1}'$ for all $n\ge1$. A counterexample to this situation can be obtained by considering  say $Y_1$ be an infinite-dimensional Banach but not Hilbert space and $Y_n=\mbox{span}[Y_1\cup\{y_1,\ldots,y_{n-1}\}]$, with $y_n\notin Y_n$ for $n\ge2$. In this case, there is no gap between any pair of $Y_n$ and $Y_{n+1}$, so it is impossible to find any $\{\tilde Y_{n_i}\}_{i\ge i_0}$ such that $\{Y_n'\}=\{Y_n\}_{n\ge1}\cup \{\tilde Y_{n_i}\}_{i\ge i_0}$ satisfies $\overline {Y_n'}\subset Y_{n+1}'$ for all $n\ge1$.
\end{remark}
\begin{remark}
There are several straightforward consequences arising  from Theorem \ref{Konyagin} which we list below.
\end{remark}
\begin{enumerate}
  \item Taking $c=1$ in Theorem \ref{Konyagin}, we obtain there exists $x\in X$ such that for $n\ge1$,
  \begin{equation}
  \label{x4}
d_n\le \rho(x,Y_n)\le 4d_n~\mbox{for all $n\ge1$}.
\end{equation}
 Note that the inequalities (\ref{x4}) improve the upper bound of $\displaystyle \frac{\rho(x,Y_n)}{d_n}$ in Theorem 1 of Konyagin \cite{Konyagin} from $8$ to $4$.  However we should point out that, Konyagin's approach in \cite{Konyagin} involves extracting a subsequence $\{d_{n_i}\}$ of $\{d_n\}$ satisfying (\ref{condition_Borodin}), as a result the element $x$ is selected such that $\rho(x,Y_n)=d_n$ for an infinite number of $n$. Our construction of $x$ in (\ref{x4}) does not necessarily satisfy this equality for infinitely many $n$, since our method involves extending $\{d_n\}$ to $\{d_n\}\cup \{K 2^{-n}\}$, and Theorem \ref{Konyagin} is applied to the \enquote{inserted sequence} $\{K2^{-n}\}$. Hence in view of our approach the obtained element $x$ satisfies $\rho(x,Y_n)=d_n$ only when the intersection $\{d_n\}_{n\ge1}\cap\{K2^{-n}\}_{n\ge i_0}$ contains an infinite number of elements.
  \item Taking $c=\displaystyle \frac{1}{4}$ in Theorem \ref{Konyagin}, we obtain existence of $x\in X$ for which
  $$
\frac{1}{4}\le \frac{\rho(x,Y_n)}{d_n}\le 1,~\mbox{for all}~n\ge1.
$$
The interval length $\displaystyle \frac{3}{4}$ makes $\big[\displaystyle \frac{1}{4},1\big]$ the \enquote{narrowest} estimating interval of $\displaystyle \frac{\rho(x,Y_n)}{d_n}$ that Theorem \ref{Konyagin} could provide.
\item Bernstein Lethargy for Fr\'{e}chet spaces is  given in \cite{Ak-Le2}, thus an improvement on Konyagin's result can be deduced. However, our proof is constructive in  explaining the relationship between the sequence of subspaces and the ``lethargic" sequence we take.
\item By using Borodin's result Theorem \ref{thm:Borodin} and under the same subspace condition on $\{Y_n\}$, we can apply the same approach as in the proof of Theorem \ref{Konyagin} to show that for any $\varepsilon>0$ arbitrarily small, there exists $x\in X$ (depending on $\varepsilon$) such that
  \begin{equation}
  \label{x5}
d_n\le \rho(x,Y_n)\le (2+\varepsilon)^2d_n~\mbox{for all $n\ge1$}.
\end{equation}
Note that (\ref{x5}) above, in fact can be obtained by replacing $\{K2^{-n}\}$ with $\{K(2+\varepsilon)^{-n}\}$ in the proof of Theorem \ref{Konyagin}. Our approximation interval of $\rho(x,Y_n)$ in (\ref{x4}) presents less deviation than  that in (\ref{x5}), due to the fact that unlike the condition (\ref{condition_Borodin}), (\ref{condition}) allows to take $d_n=K2^{-n}$, $n\ge1$.
\end{enumerate}

\noindent
{\bf Acknowledgement.} The authors are grateful to the referees and to the editor for many useful remarks which helped to make the manuscript more transparent.
\bibliographystyle{amsplain}

\end{document}